\documentclass{article}

\usepackage{amsmath}
\usepackage{amssymb}
\usepackage{amsthm}
\usepackage{graphicx}
\usepackage{bbm}

\newtheorem{theorem}{Theorem}
\newtheorem{proposition}[theorem]{Proposition}
\newtheorem{lemma}[theorem]{Lemma}

\newtheorem{example}[theorem]{Example}

\begin{document}

\title{The Positivity Set of a Recurrence Sequence}

\author{
Jason P.\ Bell
\\ Department of Mathematics
\\ University of Michigan, Ann Arbor, USA
\\ belljp@umich.edu
\and
Stefan Gerhold
 \thanks{Supported by the SFB-grant F1305 of the Austrian FWF}
 \\ RISC, J. Kepler University Linz, Austria
 \\ stefan.gerhold@risc.uni-linz.ac.at
}

\maketitle

\begin{abstract}
We consider real sequences $(f_n)$ that satisfy a linear recurrence with constant coefficients. We show
that the density of the positivity set of such a sequence always exists. In the special case
where the sequence has no positive dominating characteristic root, we establish that
the density is positive. Furthermore, we determine the values that can occur as density
of such a positivity set, both for the special case just mentioned and in general.
\end{abstract}



\section{Introduction and Main Results}

A sequence $(f_n)_{n\geq 0}$ of real numbers is called a \emph{recurrence sequence} if it satisfies
a linear recurrence
\[
  f_{n+h} = c_1 f_{n+h-1} + \dots + c_{h-1}f_{n+1} + c_h f_n, \qquad n\geq 0,
\]
with constant coefficients $c_k\in\mathbbm{R}$. (Since we are concerned with questions of positivity,
we restrict attention to real sequences.) One of the most charming and celebrated results in the
theory of recurrence sequences is the Skolem-Mahler-Lech theorem. It asserts
that the zero set
\[
  \{n\in\mathbbm{N} : f_n = 0 \}
\]
of a recurrence sequence is the union of a
finite set and finitely many arithmetic progressions. The recent comprehensive monograph by Everest et
al.~\cite{EvPoSh03}
contains references to the substantial literature devoted to this result and to related questions.
However, not much seems to be known about the positivity set
\begin{equation}\label{eq:pos set}
  \{n\in\mathbbm{N} : f_n > 0 \}.
\end{equation}
In the following section we establish that the density of the set \eqref{eq:pos set}
always exists, where the (natural) density of a set $A\subseteq \mathbbm{N}$
is defined as
\[
  \boldsymbol{d}(A):=\lim_{x\to\infty}x^{-1}\sharp\{n\leq x:n\in A\},
\]
provided that the limit exists.

\begin{theorem}\label{thm:density exists}
  Let $(f_n)$ be a recurrence sequence. Then the density of the set $\{n\in\mathbbm{N} : f_n > 0\}$ exists.
\end{theorem}

Recall that a recurrence sequence can be written
as a generalized power sum
\begin{equation}\label{eq:power sum}
  f_n = \sum_{k=0}^m P_k(n) \gamma_k^n, \qquad n\geq 0,
\end{equation}
with non-zero polynomials $P_k\in\mathbbm{C}[n]$ and \emph{roots} $\gamma_k\in\mathbbm{C}$
that are roots of the \emph{characteristic polynomial}
\[
  z^h - c_1 z^{h-1} - \dots - c_{h-1} z - c_h
\]
of the recurrence. The roots of largest modulus are called
\emph{dominating roots} of $(f_n)$.

It does not come as a surprise that recurrence sequences with no positive dominating root have oscillating behaviour.
Indeed, in section~\ref{se:no pos root} we prove that for such a sequence $(f_n)$
the densities of \eqref{eq:pos set} and the negativity set
\begin{equation}\label{eq:neg set}
  \{n\in\mathbbm{N} : f_n < 0 \}
\end{equation}
are always positive. This generalizes the following known result~\cite{BuWe81,Ge05}:
If a recurrence sequence has at most four dominating roots, and none of them
is real positive, then the sets \eqref{eq:pos set} and \eqref{eq:neg set}
both have infinitely many elements.

\begin{theorem}\label{thm:no pos root}
  Let $(f_n)$ be a nonzero recurrence sequence with no positive dominating characteristic root. Then the
  sets $\{n\in\mathbbm{N}:f_n>0\}$ and $\{n\in\mathbbm{N}:f_n<0\}$ have positive density.
\end{theorem}

In section~\ref{se:poss dens} we investigate which numbers actually occur
as density of the positivity set of some recurrence sequence. It turns out
that all possible values occur, both for sequences with no dominating positive root
and in general.

Finally (section~\ref{se: weak SML}), we return to the Skolem-Mahler-Lech theorem.
Our approach yields the following weak version: The density of the zero set exists
and is a rational number.

The conclusion hints at algorithmic aspects of the positivity of recurrence sequences.

\section{The Density of the Positivity Set}

Notation: We write $f_n\equiv 0$ if $f_n=0$ for all $n\geq 0$. The Lebesgue measure
of a set $B\subset\mathbbm{R}^m$ is denoted by $\boldsymbol{\lambda}(B)$.

The goal of this section is to prove Theorem~\ref{thm:density exists}.
Dividing $f_n$ by $n^D|\gamma_1|^n$, where $\gamma_1$ is a dominating root of $f_n$ and $D$ is
the maximal degree of the $P_k$ with $|\gamma_k|=|\gamma_1|$, we obtain from \eqref{eq:power sum}
\[
  n^{-D}|\gamma_1|^{-n} f_n = \sum_{i=1}^d a_i \cos(2\pi \theta_i n + \beta_i) + v - r_n,
\]
where $r_n=\mathrm{O}(1/n)$ is a recurrence sequence,
$\theta_1,\dots,\theta_d$ are in $]0,1[$ and $a_i,\beta_i,v\in\mathbbm{R}$.
From now on we will assume w.l.o.g.\ $D=0$ and $|\gamma_1|=1$. As a first step
we get rid of any integer relations that the $\theta_i$'s might satisfy.

\begin{lemma}\label{le:module basis}
  Let $\theta_1,\dots,\theta_d$ be real numbers. Then there is a basis $\{\tau_1,\dots,\tau_{m+1}\}$
  of the $\mathbbm{Z}$-module
  \[
    M = \mathbbm{Z} + \mathbbm{Z}\theta_1 + \dots \mathbbm{Z}\theta_d
  \]
  such that $1/\tau_{m+1}$ is a positive integer and $1,\tau_1,\dots,\tau_m$ are linearly independent over $\mathbbm{Q}$.
\end{lemma}
\begin{proof}
  $M$ is finitely generated and torsion free, hence it is free~\cite[Theorem~III.7.3]{La02}.
  Let $\{\alpha_1,\dots,\alpha_{m+1}\}$ be a basis.
  Since $1\in M$, there are integers $e_1,\dots,e_{m+1}$ such that
  \[
    e_1 \alpha_1 + \dots + e_{m+1} \alpha_{m+1} = 1.
  \]
  We complete $(e_1/g,\dots,e_{m+1}/g)$, where $g:=\gcd(e_1,\dots,e_{m+1})$,
  to a unimodular integer matrix ${\bf C}$ with last row $(e_1/g,\dots,e_{m+1}/g)$ ~\cite[\S XXI.3]{La02}.
  Then
  \[
    (\tau_1,\dots,\tau_{m+1})^T := {\bf C}(\alpha_1,\dots,\alpha_{m+1})^T
  \]
  yields a basis of $M$ with $\tau_{m+1}=1/g\in\mathbbm{Q}$. Now suppose
  \[
    u_1 \tau_1 + \dots + u_m \tau_m = u
  \]
  for integers $u_1,\dots,u_m,u$. Since $u$ has also the representation
  \[
    u g \tau_{m+1} = u,
  \]
  it follows $u_1=\dots=u_m=u=0$.
\end{proof}

Take $\tau_1,\dots,\tau_{m+1}$ as in Lemma~\ref{le:module basis}, with $\tau_{m+1}=1/g$.
Roughly speaking, we have put all integer relations among the $\theta_i$ into the rational basis element $\tau_{m+1}$.
There are integers $b_{ij}$
with
\[
  \theta_i = \sum_{j=1}^{m+1} b_{ij}\tau_j.
\]
Now we split the sequence $(f_n)$ into the subsequences $(f_{gn+k})_{n\geq 0}$
for $0\leq k<g$. We have
\[
  f_{gn+k} = G_n - s_n,
\]
where $s_n:=r_{gn+k}$ and $G_n$ is the dominant part.
Defining the integer matrix
\[
  {\bf B} := (g b_{ij})_{\substack{1\leq i\leq d \\ 1\leq j\leq m}} \in \mathbbm{Z}^{d\times m}
\]
and the real vector ${\bf c} = (c_1,\dots,c_d)$ with
\[
  c_i := 2\pi k\sum_{j=1}^{m+1}b_{ij}\tau_j+\beta_i, \qquad 1\leq i\leq d,
\]
it can be written as (cos is applied component wise)
\[
  G_n = {\bf a}^T \cos(2\pi n {\bf B}\boldsymbol{\tau} + {\bf c}) + v.
\]
We show that the density of $\{n\in\mathbbm{N} : f_{gn+k} > 0\}$ exists for each $k$.
Since $s_n$ is a recurrence sequence with fewer characteristic roots than $f_n$,
we may assume inductively that $\boldsymbol{d}(\{n\in\mathbbm{N}:s_n<0\})$ exists.
Thus, if $G_n$ is the zero sequence, we are done. Now let $k$ be such that $G_n$ is not the zero sequence.
It is plain that $G_n=H(n\boldsymbol{\tau})$, where 
\[
  H({\bf t}) := {\bf a}^T \cos(2\pi {\bf B}{\bf t} + {\bf c}) + v, \quad {\bf t}\in [0,1]^m.
\]
The following theorem shows that the function $H$ can be used to evaluate the density
of the positivity set of $G_n$, which equals, as we will see below, that of the set $\{n:f_{gn+k}>0\}$.
\begin{theorem}[Kronecker-Weyl]\label{thm:kronecker}
  Let $\tau_1,\dots,\tau_m$ be real numbers such that $1,\tau_1,\dots,\tau_m$ are linearly independent over $\mathbbm{Q}$.
  Then for every Jordan measurable set $A\subseteq [0,1]^m$ we have
  \[
    \boldsymbol{d}(\{n\in\mathbbm{N} : n\boldsymbol{\tau}\bmod 1 \in A\}) = \boldsymbol{\lambda}(A).
  \]
\end{theorem}
\begin{proof}
  We refer to Cassels~\cite[Theorems~IV.I and IV.II]{Ca57}.
\end{proof}
We define
\begin{equation}\label{eq:sets G_n}
  L_\varepsilon := \{ n\in\mathbbm{N} : G_n \geq \varepsilon \} \qquad \text{and} \qquad
  S_\varepsilon := \{ n \in\mathbbm{N}: |G_n| < \varepsilon \}.
\end{equation}
The corresponding sets for the function $H$ are defined as
\[
  \widetilde{L}_\varepsilon := \{ {\bf t}\in [0,1]^{m} : H({\bf t}) \geq \varepsilon \} \qquad \text{and} \qquad
  \widetilde{S}_\varepsilon := \{ {\bf t}\in [0,1]^{m} : |H({\bf t})| < \varepsilon \}.
\]
Since for all $\varepsilon \geq 0$
\[
  L_\varepsilon
  = \{ n \in \mathbbm{N} : n\boldsymbol{\tau} \bmod 1 \in \widetilde{L}_\varepsilon \},
\]
we have
$\boldsymbol{d}(L_\varepsilon) = \boldsymbol{\lambda}(\widetilde{L}_\varepsilon)$
for all $\varepsilon \geq 0$ by Theorem~\ref{thm:kronecker}. Similarly,
\begin{equation}\label{eq:S_eps}
  \boldsymbol{d}(S_\varepsilon) = \boldsymbol{\lambda}(\widetilde{S}_\varepsilon),\quad \varepsilon > 0.
\end{equation}
Note that the boundary of the bounded set $\widetilde{S}_\varepsilon$ (respectively $\widetilde{L}_\varepsilon$) is a Lebesgue null set
(as seen by applying the following lemma with $F({\bf t})=H({\bf t})-\varepsilon$),
hence $\widetilde{S}_\varepsilon$ and $\widetilde{L}_\varepsilon$ are Jordan measurable, and Theorem~\ref{thm:kronecker}
is indeed applicable. Lemma~\ref{le:zero set} seems to be 
known~\cite{JiSiBe05}, but we could not find a complete proof in the literature.

\begin{lemma}\label{le:zero set}
  Let $F:\mathbbm{R}^m \rightarrow \mathbbm{R}$ be a real analytic function.
  Then the zero set of $F$ has Lebesgue measure zero, unless $F$ vanishes identically.
\end{lemma}

The proof of Lemma~\ref{le:zero set} is postponed to the end of this section.
Since $G_n$ is not the zero
sequence, the function $H$ does not vanish identically on $[0,1]^m$.
By the Lebesgue dominated convergence theorem and Lemma~\ref{le:zero set} we thus find
\[
  \lim_{\varepsilon\to 0} \boldsymbol{\lambda}(\widetilde{S}_\varepsilon) = 0 \qquad \text{and} \qquad
  \lim_{\varepsilon\to 0} \boldsymbol{\lambda}(\widetilde{L}_\varepsilon) = \boldsymbol{\lambda}(\widetilde{L}_0).
\]
This yields $\boldsymbol{d}(\{n\in\mathbbm{N} : G_n > s_n \}) = \boldsymbol{\lambda}(\widetilde{L}_0)$ by the
following lemma, which completes the proof of Theorem~\ref{thm:density exists}.

\begin{lemma}\label{le:density}
  Let $G_n$ and $s_n$ be real sequences with $s_n=\mathrm{o}(1)$ and let $L_\varepsilon$, $S_\varepsilon$  be as in \eqref{eq:sets G_n}.
  Suppose that $\boldsymbol{d}(L_\varepsilon)$ and $\boldsymbol{d}(S_\varepsilon)$
  exist for all $\varepsilon\geq 0$, and that
  \[
    \lim_{\varepsilon\to 0} \boldsymbol{d}(L_\varepsilon) = \boldsymbol{d}(L_0) \qquad and \qquad
    \lim_{\varepsilon\to 0} \boldsymbol{d}(S_\varepsilon) = 0.
  \]
  Then
  \[
    \boldsymbol{d}(\{ n\in\mathbbm{N} : G_n > s_n\}) = \boldsymbol{d}(L_0).
  \]
\end{lemma}
\begin{proof}
  For any set $A\subseteq \mathbbm{N}$ we write $A(x):=\{n\leq x: n \in A\}$. 
  Define
  \[
    P:=\{n\in\mathbbm{N} : G_n > s_n \}.
  \]
  Let $\varepsilon>0$ be arbitrary. Take $n_0$ such
  that $|s_n|<\varepsilon$ for $n > n_0$. It follows

  \begin{align*}
    \sharp P(x) &= \sharp\{n \leq n_0 : G_n>s_n \} + \sharp\{n_0 < n\leq x : G_n \geq \varepsilon \} \\
    &\phantom{=} +\sharp\{n_0 < n\leq x : s_n<G_n < \varepsilon \},
  \end{align*}
  hence
  \[
    |\sharp P(x) - \sharp L_\varepsilon(x) | \leq \sharp S_\varepsilon(x)+\mathrm{o}(x)
  \]
  as $x\to \infty$. Thus we have
  \begin{align*}
    |x^{-1}\sharp P(x) - \boldsymbol{d}(L_0)| &\leq |x^{-1}\sharp P(x) - x^{-1}\sharp L_\varepsilon(x)| 
      + |x^{-1}\sharp L_\varepsilon(x) - \boldsymbol{d}(L_0)| \\
    & \leq x^{-1}\sharp S_\varepsilon(x) + |x^{-1}\sharp L_\varepsilon(x) - \boldsymbol{d}(L_0)| + \mathrm{o}(1).
  \end{align*}
  The right hand side goes to
  \begin{equation*}
    \boldsymbol{d}(S_\varepsilon)+|\boldsymbol{d}(L_\varepsilon) - \boldsymbol{d}(L_0)|
  \end{equation*}
  as $x\to\infty$.
  By assumption, this can be made arbitrarily small, which implies $\boldsymbol{d}(P) = \boldsymbol{d}(L_0)$.
\end{proof}

\begin{proof}[Proof of Lemma~\ref{le:zero set}]
  For $m=1$ this is clear, since then the zero set is countable. Now assume that we have established the result
  for $1,\dots,m-1$. Put
  \[
    V := \{ (t_2,\dots,t_m) \in \mathbbm{R}^{m-1} : F(\cdot,t_2,\dots,t_m) \text{ vanishes identically }\}.
  \]
  Take a real number $s$ such that $F(s,\cdot,\dots,\cdot)$ is not identically zero. Clearly, $F(s,t_2,\dots,t_m)=0$
  for all $(t_2,\dots,t_m)\in V$. By the induction hypothesis, this implies $\boldsymbol{\lambda}(V)=0$. Note that $V$ is closed, hence measurable.
  Since $F$ is real analytic in the first argument, we have
  \[
    \int_\mathbbm{R} \chi_Z(t_1,\dots,t_m) \mathrm{d}\boldsymbol{\lambda}(t_1) = 0
  \]
  for all $(t_2,\dots,t_m)\notin V$, where $\chi_Z$ is the characteristic function of the zero set
  \[
    Z: = \{ (t_1,\dots,t_m) \in \mathbbm{R}^m : F(t_1,\dots,t_m)=0 \}.
  \]
  Since $V$ has measure zero, this implies
  \[
    \int_\mathbbm{R}\dots\int_\mathbbm{R} \chi_Z(t_1,\dots,t_m) \mathrm{d}\boldsymbol{\lambda}(t_1) \cdots \mathrm{d}\boldsymbol{\lambda}(t_m) = 0.
  \]
  This argument works for any order of integration, hence we obtain $\int_{\mathbbm{R}^m}\chi_Z=0$ by Tonelli's theorem.
\end{proof}

\section{Sequences with no Positive Dominating Root}\label{se:no pos root}

In this section we prove Theorem~\ref{thm:no pos root}. We begin by settling the
special cases where the $\theta_i$ are all irrational or all rational, and then
put them together.

\begin{lemma}\label{le:irr}
  Let $\theta_1,\dots,\theta_d$ be irrational numbers,
  and let $a_i,\beta_i$ be real numbers such that the sequence
  \[
    u_n = \sum_{i=1}^d a_i \cos(2\pi \theta_i n + \beta_i)
  \]
  is not identically zero. Let further $r_n$ be a recurrence sequence with $r_n=\mathrm{o}(1)$.
  Then the set $\{n \in \mathbbm{N} : u_n > r_n\}$ has positive density.
\end{lemma}
\begin{proof}
  Proceeding as in the proof of Theorem~\ref{thm:density exists}, we can write
  \[
    G_n := u_{gn+k} = {\bf a}^T \cos(2\pi n {\bf B}\boldsymbol{\tau} + {\bf c}),
  \]
  where ${\bf B}$ is an integer matrix no row of which is zero, ${\bf c}$ is a real vector and $1,\tau_1,\dots,\tau_m$ are linearly
  independent over $\mathbbm{Q}$.
  If $k$ is such that $G_n=u_{gn+k}\equiv 0$, then the density of $\{n\in\mathbbm{N}:G_n > s_n\}$, where
  $s_n=r_{qn+k}$, exists by Theorem~\ref{thm:density exists}, but may be zero.
  Now choose a $k_0$ such that the corresponding
  sequence $G_n=u_{gn+k_0}$ is not the zero sequence.
  We have $G_n=H(n\boldsymbol{\tau})$, where
  \[
    H({\bf t}) := {\bf a}^T \cos(2\pi {\bf B}{\bf t} + {\bf c}).
  \]
  Moreover, with the notation of the proof of Theorem~\ref{thm:density exists}, we have
  \[
    \boldsymbol{d}(\{n\in\mathbbm{N} : G_n > s_n \}) = \boldsymbol{\lambda}(\widetilde{L}_0).
  \]
  The function $H$ is not identically zero on $[0,1]^m$. But
  \begin{equation}\label{eq:int H}
    \int_0^1\dots\int_0^1 H(t_1,\dots,t_m) \mathrm{d}t_1 \cdots \mathrm{d}t_m = 0,
  \end{equation}
  because no row of ${\bf B}$ is the zero vector.
  Hence $H$ has a positive value on $[0,1]^m$, and since it is continuous, we have $\boldsymbol{\lambda}(\widetilde{L}_0)>0$.
\end{proof}

Observe that the integral in \eqref{eq:int H} need not vanish if ${\bf B}$ has a zero row, which can
only happen if the $\theta_i$ corresponding to this row is a rational number. This is the reason
why we consider rational $\theta_i$'s separately.

\begin{lemma}\label{le:rat}
  Let $\theta_1,\dots,\theta_d$ be rational numbers in $]0,1[$,
  and let $a_i,\beta_i$ be real numbers such that the purely periodic sequence
  \[
    u_n = \sum_{i=1}^d a_i \cos(2\pi \theta_i n + \beta_i)
  \]
  is not identically zero. Then $u_n$ has a positive and a negative value.
\end{lemma}
\begin{proof}
  By the identity
  \[
    \sum_{k=0}^{q-1}\cos{\tfrac{2\pi k p}{q}} + \mathrm{i} \sum_{k=0}^{q-1}\sin{\tfrac{2\pi k p}{q}}
    = \sum_{k=0}^{q-1}\mathrm{e}^{2\pi\mathrm{i} k p/q} = 0,
  \]
  valid for integers $0<p<q$, and the addition formula of cos we obtain
  \[
    u_0+\dots+u_{q-1}=0,
  \]
  where $q$ is a common denominator of $\theta_1,\dots,\theta_d$.
\end{proof}

\begin{proof}[Proof of Theorem~\ref{thm:no pos root}]
  It suffices to consider the positivity set.
  We may write
  \[
    f_n = u_n + v_n - r_n,
  \]
  where $r_n=\mathrm{o}(1)$ is a recurrence sequence,
  \begin{align*}
    u_n &= \sum_{i=1}^d a_i \cos(2\pi \theta_i n + \beta_i), \\
    v_n &= \sum_{i=d+1}^e a_i \cos(2\pi \theta_i n + \beta_i),
  \end{align*}
  $\theta_1,\dots,\theta_d$ are irrational, $\theta_{d+1},\dots,\theta_e$ are rational numbers in $]0,1[$
  with common denominator $q>0$ and $u_n+v_n\not\equiv 0$.
  If $v_n\equiv 0$, then the result follows from Lemma~\ref{le:irr}. Now suppose $v_n\not\equiv 0$.
  Then for each $k$ the density of the set $\{n\in\mathbbm{N}:f_{qn+k}>0\}$
  exists by Theorem~\ref{thm:density exists}. By Lemma~\ref{le:rat} there is $k_0$ such that $v_{qn+k_0}=v>0$.
  It suffices to show that the set $\{n\in\mathbbm{N}:f_{qn+k_0}>0\}$ has positive density.
  This is clear if $u_{qn+k_0}\equiv 0$. Otherwise, notice that
  \[
    \{n\in\mathbbm{N}:f_{qn+k_0}>0\} \supseteq \{n\in\mathbbm{N}:u_{qn+k_0} > r_{qn+k_0}\},
  \]
  and the latter set has positive density by Lemma~\ref{le:irr}.
\end{proof}

\section{The Possible Values of the Density}\label{se:poss dens}

In this section we investigate which values from $[0,1]$ occur as density of some
recurrence sequence. In its basic form, the question is readily answered:

\begin{example}\label{ex:poss dens}
  Let $w$ be a real number and define
  \[
    f_n := \sin(2\pi n \sqrt{2}) - w.
  \]
  Then, by Theorem~\ref{thm:kronecker},
  \begin{align*}
    \boldsymbol{d}(\{n\in\mathbbm{N}:f_n>0\}) &= \boldsymbol{\lambda}(\{t \in [0,1] : \sin(2\pi t) > w \}) \\ &=
    \begin{cases}
      1 & w\leq -1 \\
      \tfrac{1}{2} - \tfrac{1}{\pi} \arcsin{w} & -1 \leq w \leq 1 \\
      0 & w\geq 1
    \end{cases}.
  \end{align*}
  Since the range of $\arcsin$ is $[-\tfrac{\pi}{2},\tfrac{\pi}{2}]$,
  for every $\kappa \in [0,1]$ this yields a recurrence sequence $f_n$ such that
  \[
    \boldsymbol{d}(\{n\in\mathbbm{N}:f_n>0\})=\kappa.
  \]
\end{example}

The following proposition generalizes this example. Note that the density of the
zero set of a recurrence sequence is always a rational number by
the Skolem-Mahler-Lech theorem.

\begin{proposition}
  Let $\kappa$ be a real number and $r$ be a rational number with $0\leq \kappa,r \leq 1$
  and $\kappa + r \leq 1$. Then there is a recurrence sequence $(f_n)$ such that
  \[
    \boldsymbol{d}(\{n\in\mathbbm{N}:f_n>0\})=\kappa \qquad \text{and} \qquad
    \boldsymbol{d}(\{n\in\mathbbm{N}:f_n=0\})=r.
  \]
\end{proposition}
\begin{proof}
  Suppose that $r=p/q$ for positive integers $p$ and $q$.
  As seen in Example~\ref{ex:poss dens}, there is a recurrence sequence $(g_n)$ such that the density
  of the zero set of $(g_n)$ is zero and the density of its positivity set is $\kappa/(1-r)$ (The case $r=1$ is trivial).
  The interlacing sequence
  \[
    f_{bn+k} :=
    \begin{cases}
      0 & 0\leq k<p \\
      g_n & p\leq k <q 
    \end{cases}
  \]
  is a recurrence sequence~\cite[section~4.1]{EvPoSh03}. Clearly, the density of its zero set is $r$,
  and the density of its positivity set is
  \[
    \boldsymbol{d}(\{n\in\mathbbm{N}:f_n>0\}) = \frac{q-p}{q} \times \frac{\kappa}{1-r} = \kappa,
  \]
  as required.
\end{proof}

If we restrict attention to sequences without dominating real positive roots, then Theorem~\ref{thm:no pos root}
tells us that the density of the positivity set can be neither zero nor one. Still, all
values in between occur.

\begin{theorem}
  Let $\kappa \in ]0,1[$. Then there is a recurrence sequence $(f_n)$ with no positive dominating characteristic root and
  $\boldsymbol{d}(\{n\in\mathbbm{N}:f_n>0\})=\kappa$.
\end{theorem}
\begin{proof}
  Let $\varepsilon>0$ be arbitrary.
  We define a function $H$ on $[0,\tfrac{1}{2}]$ by
  \[
    H(t) :=
    \begin{cases}
      \tfrac{(\varepsilon-1)^2}{\varepsilon} \left( 1-\tfrac{2t}{\varepsilon} \right) & 0\leq t\leq \tfrac{\varepsilon}{2} \\
      \varepsilon - 2t & \tfrac{\varepsilon}{2} \leq t\leq \tfrac{1}{2}
    \end{cases}
  \]
  and extend it to an even, $1$-periodic function $H$ on $\mathbbm{R}$ (see Figure~\ref{fig:H}). It is continuous and satisfies
  \[
    \int_0^1 H(t)\mathrm{d}t=0 \qquad \text{and} \qquad
    \boldsymbol{\lambda}(\{ t\in[0,1] : H(t) > 0 \}) = \varepsilon.
  \]
  
  \begin{figure}[h]
    \centering
    \includegraphics[scale=.7]{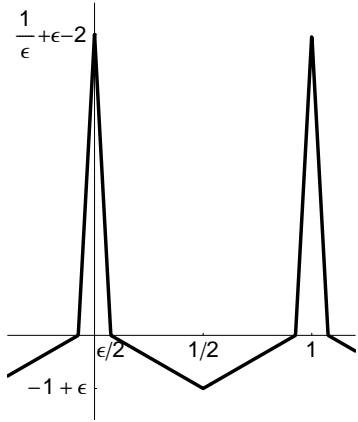}
    \caption{The function $H$}
    \label{fig:H}
  \end{figure}

  Expanding $H$ into a Fourier series, we find that there are real $a_j$ such that $H$ is the pointwise limit of
  \[
    H_m(t) := \sum_{j=1}^m a_j \cos(2\pi j t)
  \]
  as $m\to\infty$. Since the zero set of $H$ is a null set, the Lebesgue dominated convergence theorem yields
  \[
    \lim_{m\to\infty} \boldsymbol{\lambda}(\{ t\in[0,1] : H_m(t) > 0 \}) = \varepsilon.
  \]
  We fix an $m$ such that
  \[
    \boldsymbol{\lambda}(\{ t\in[0,1] : H_m(t) > 0 \}) \leq 2\varepsilon.
  \]
  The function
  \[
    \phi(A_1,\dots,A_m) := \boldsymbol{\lambda}\Bigl(\Bigl\{ t\in[0,1] : \sum_{j=1}^m A_j \cos(2\pi j t) > 0 \Bigr\}\Bigr)
  \]
  is continuous on $\mathbbm{R}^m\backslash\{{\bf 0}\}$. To see this, observe that $\phi$ is continuous
  at all points $(A_1,\dots,A_m)$ for which
  $\sum_{j=1}^m A_j \cos(2\pi j t)$ is not identically zero and appeal to the uniqueness of the Fourier expansion.
  Since $\phi(1,0,\dots,0)=\tfrac{1}{2}$
  and $\phi(a_1,\dots,a_m) \leq 2\varepsilon$, the function $\phi$
  assumes every value from $[2\varepsilon,\tfrac{1}{2}]$ by the intermediate value theorem.

  Hence the positivity sets of the sequences
  \[
    f_n := \sum_{j=1}^m A_j \cos(2\pi j n \sqrt{2})
  \]
  assume
  all densities from $[2\varepsilon,\tfrac{1}{2}]$ for appropriate choices of $(A_1,\dots,A_m)$ by Theorem~\ref{thm:kronecker}.
  Repeating the whole argument with $-H$ instead of $H$ yields the desired result for
  $\kappa \in [\tfrac{1}{2},1-2\varepsilon]$. Since $\varepsilon$ was arbitrary, the theorem is proved.
\end{proof}

\section{A Weak Version of Skolem-Mahler-Lech}\label{se: weak SML}

Without using the Skolem-Mahler-Lech theorem, it follows from Theorem~\ref{thm:density exists}
that the density of the zero set of a recurrence sequence $(f_n)$ exists. We can show a bit
more with our approach. Recall, however, that we only deal with real sequences, whereas
the Skolem-Mahler-Lech theorem holds for any field of characteristic zero.

\begin{proposition}
  The density of the zero set of a (real) recurrence sequence $(f_n)$ is a rational number.
\end{proposition}
\begin{proof}
  Let $k$ be a natural number, and let $g$, $G_n$ and $s_n$ be as in the proof of Theorem~\ref{thm:density exists}.
  If $k$ is such that $G_n\equiv 0$, then the density of the zero set of $f_{gn+k}$ is rational,
  since we may assume inductively that the density of $\{n:s_n=0\}$ is rational.

  Now suppose $G_n\not\equiv 0$. The zero set of $f_{gn+k}$ can be partitioned as
  \[
    \{n\in\mathbbm{N}:G_n=s_n\} = \{n : G_n=s_n, |G_n|<\varepsilon \} \cup
      \{n : G_n=s_n, |G_n|\geq\varepsilon \},
  \]
  where $\varepsilon \geq 0$ is arbitrary.
  The latter set is finite, and the first one is contained in $S_\varepsilon$, defined
  in \eqref{eq:sets G_n}. Hence
  \[
    \boldsymbol{d}(\{n\in\mathbbm{N}:G_n=s_n\}) \leq \boldsymbol{d}(S_\varepsilon)
  \]
  for all $\varepsilon \geq 0$. But we know that $\lim_{\varepsilon\to 0}\boldsymbol{d}(S_\varepsilon) = 0$
  from the proof of Theorem~\ref{thm:density exists}, which yields
  \[
    \boldsymbol{d}(\{n\in\mathbbm{N}:G_n=s_n\}) = 0.
  \]
  Thus, the zero sets of all subsequences $(f_{gn+k})_{n\geq 0}$, $0\leq k<g$, have rational density,
  which proves the desired result.
\end{proof}

\section{Conclusion}

There is no algorithm known for deciding whether $f_n>0$ for all $n$, nor has the problem been shown to be undecidable.
When we are talking about algorithmics, it is natural to assume that the recurrence coefficients and
the initial values are rational numbers.
In this case Gourdon and Salvy~\cite{GoSa:96} have proposed an efficient method for ordering the characteristic roots
w.r.t.\ to their modulus. Thus, the dominating characteristic roots can be identified algorithmically.
If none of them is real positive, then we know that the sequence oscillates by Theorem~\ref{thm:no pos root}.
On the other hand, sequences where a positive dominating root is accompanied by complex dominating roots
seem to pose difficult Diophantine problems.
For instance, we do not know if the sequence
\begin{equation}\label{eq:difficult}
  f_n := \cos(2\pi\theta n) + 1 + \left(-\tfrac{1}{2}\right)^n
\end{equation}
is positive for $\theta=\sqrt{2}$, say. It can be shown, however, that the set
of $\theta$'s for which the corresponding sequence $(f_n)$ (defined by~\eqref{eq:difficult}) is positive has measure zero~\cite[Theorem~7.2]{Ge05}.

\bibliographystyle{siam}
\bibliography{/home/sgerhold/documents/mypapers/gerhold}

\end{document}